\documentclass[11pt]{article}

\usepackage{amsmath}          % various ams stuf
\usepackage{amssymb}          % dito
\usepackage{latexsym}         % dito
\usepackage{amsxtra}          % dito
\usepackage[latin1]{inputenc} % smart input of special characters
\usepackage{eucal}            % nicer caligraphic fonts
\usepackage{bbm}              % nicer bbm fonts
\usepackage{exscale}          % larger summation signs
\usepackage{theorem}          % nicer theorems
\usepackage[curve]{xypic}     % the cool diagrams need this with pdf
\usepackage{diagxy}           % cool diagrams
\usepackage[sort]{cite}       % nicer citations
\usepackage{mathrsfs}         % nicer caligraphic fonts
\usepackage{stmaryrd}         % more symbols & brackets
\title{Morita Theory in Deformation Quantization\\[0.6cm]
  \begin{normalsize}
      Dedicated to the memory of Nikolai Neumaier
  \end{normalsize}
}
\author{\textbf{Stefan Waldmann}\thanks{\texttt{Stefan.Waldmann@physik.uni-freiburg.de}}
  \\[0.2cm]
  \begin{minipage}{8cm}
      \begin{center}
          \begin{small}
              Fakult{\"a}t f{\"u}r Mathematik und Physik \\
              Albert-Ludwigs-Universit{\"a}t Freiburg \\
              Physikalisches Institut \\
              Hermann Herder Strasse 3 \\
              D 79104 Freiburg \\
              Germany
          \end{small}
    \end{center}
\end{minipage}
}
\date{December 2010}

\renewcommand{\mathbb}[1]{\mathbbm{#1}}          % use nicer bbm fonts
\allowdisplaybreaks                              % good for long formulas
\renewcommand{\arraystretch}{1.2}                % in formulas and tables

\theoremheaderfont{\normalfont\bfseries}
\theorembodyfont{\itshape}
\newtheorem{lemma}{Lemma}[section]
\newtheorem{theorem}[lemma]{Theorem}
\newtheorem{definition}[lemma]{Definition}

\theorembodyfont{\rmfamily}

\newtheorem{example}[lemma]{Example}
\newtheorem{remark}[lemma]{Remark}

\newcommand{\Lie}        {\operatorname{\mathscr{L}\!}}    
\newcommand{\cc}[1]      {\overline{{#1}}}              
\newcommand{\id}         {\operatorname{\mathsf{id}}}   
\newcommand{\End}        {\operatorname{\mathsf{End}}}   
 
\newcommand{\Unit}       {\mathbb{1}}                  
\newcommand{\cl}         {\mathrm{cl}}                    
\newcommand{\I}          {\mathrm{i}}
\newcommand{\E}          {\mathrm{e}}
\newcommand{\D}          {\operatorname{\mathrm{d}}} 
\newcommand{\Anti}       {\Lambda}
\newcommand{\Schouten}[1]{\left\llbracket{#1}\right\rrbracket}
\newcommand{\lie}[1]     {\mathfrak{#1}}
\newcommand{\ring}[1]    {\mathsf{#1}}
\newcommand{\deform}[1]  {\boldsymbol{#1}}

\newcommand{\Diffop}     {\operatorname{\mathrm{DiffOp}}}

\newcommand{\Bimod}[5] {\sideset{^{\scriptscriptstyle{#1}}_{\scriptscriptstyle{#2}}}{^{\scriptscriptstyle{#4}}_{\scriptscriptstyle{#5}}}{\operatorname{#3}}}

\newcommand{\EA}   {\Bimod{}{}{\mathcal{E}}{}{\mathcal{A}}}
\newcommand{\BEpA}  {\Bimod{}{\mathcal{B}}{\mathcal{E}}{\prime}{\mathcal{A}}}
\newcommand{\BEA}  {\Bimod{}{\mathcal{B}}{\mathcal{E}}{}{\mathcal{A}}}
\newcommand{\CFB}  {\Bimod{}{\mathcal{C}}{\mathcal{F}}{}{\mathcal{B}}}
\newcommand{\AAA}  {\Bimod{}{\mathcal{A}}{\mathcal{A}}{}{\mathcal{A}}}
\newcommand{\BBB}  {\Bimod{}{\mathcal{B}}{\mathcal{B}}{}{\mathcal{B}}}

\newcommand{\IP}[4]{{\,}_{\scriptscriptstyle{#2}\!\!}\left\langle{{#1}}\right\rangle^{\scriptscriptstyle{#3}}_{\scriptscriptstyle{#4}}}

\newcommand{\SPA}[1]     {\IP{{#1}}{}{}{\mathcal{A}}}
\newcommand{\SPEA}[1]    {\IP{{#1}}{}{\mathcal{E}}{\mathcal{A}}}
\newcommand{\SPFB}[1]    {\IP{{#1}}{}{\mathcal{F}}{\mathcal{B}}}
\newcommand{\SPFEA}[1]   {\IP{{#1}}{}{\mathcal{F} \tensor \mathcal{E}}{\mathcal{A}}}

\newcommand{\tensor}[1][{}] {\mathbin{\otimes_{\scriptscriptstyle{#1}}}}
\newcommand{\itensor}[1][{}]{\mathbin{\widehat{\otimes}_{\scriptscriptstyle{#1}}}}

\newcommand{\Pic}      {\operatorname{\mathsf{Pic}}}
\newcommand{\PicBI}    {\operatorname{\underline{\mathsf{Pic}}}}
\newcommand{\PicH}      {\sideset{}{_{H}}{\operatorname{\mathsf{Pic}}}}
\newcommand{\StrPic}   {\sideset{}{^{\mathrm{str}}}{\operatorname{\mathsf{Pic}}}}
\newcommand{\starPic}  {\sideset{}{^*}{\operatorname{\mathsf{Pic}}}}
\newcommand{\StrPicH}  {\sideset{}{^{\mathrm{str}}_{H}}{\operatorname{\mathsf{Pic}}}}
\newcommand{\starPicH} {\sideset{}{^*_{H}}{\operatorname{\mathsf{Pic}}}}
\newcommand{\BiMod}       {\operatorname{\mathsf{Bimod}}}
\newcommand{\starBiMod}   {\sideset{}{^{*}}{\operatorname{\mathsf{Bimod}}}}
\newcommand{\StrBiMod}    {\sideset{}{^{\mathrm{str}}}{\operatorname{\mathsf{Bimod}}}}
\newcommand{\BiModH}       {\sideset{}{_{H}}{\operatorname{\mathsf{Bimod}}}}
\newcommand{\starBiModH}   {\sideset{}{^{*}_{H}}{\operatorname{\mathsf{Bimod}}}}
\newcommand{\StrBiModH}    {\sideset{}{^{\mathrm{str}}_{H}}{\operatorname{\mathsf{Bimod}}}}
\newcommand{\BiModBI}     {\operatorname{\underline{\mathsf{Bimod}}}}
\newcommand{\starBiModBI} {\sideset{}{^{*}}{\operatorname{\underline{\mathsf{Bimod}}}}}
\newcommand{\StrBiModBI}  {\sideset{}{^{\mathrm{str}}}{\operatorname{\underline{\mathsf{Bimod}}}}}
\newcommand{\BiModBIH}     {\sideset{}{^{*}_{H}}{\operatorname{\underline{\mathsf{Bimod}}}}}
\newcommand{\starBiModBIH} {\sideset{}{^{*}_{H}}{\operatorname{\underline{\mathsf{Bimod}}}}}
\newcommand{\StrBiModBIH}  {\sideset{}{^{\mathrm{str}}_{H}}{\operatorname{\underline{\mathsf{Bimod}}}}}

\newcommand{\FormalPoiss}    {\operatorname{\underline{FPoiss}}}
\newcommand{\Deformation}    {\operatorname{\underline{Def}}}
\newcommand{\FDiffeo}    {\operatorname{FDiffeo}}
\newcommand{\Def}    {\operatorname{Def}}
\newcommand{\FPoiss}     {\operatorname{FPoiss}}
\newcommand{\HdR}         {\mathrm{H}_{\scriptscriptstyle \mathrm{dR}}}
\newcommand{\acts}       {\mathbin{\triangleright}}
\newcommand{\sweedler}[1] {{\scriptscriptstyle{(#1)}}}

\begin{document}

\maketitle

\begin{abstract}
    Various aspects of Morita theory of deformed algebras and in
    particular of star product algebras on general Poisson manifolds
    are discussed. We relate the three flavours ring-theoretic Morita
    equivalence, $^*$-Morita equivalence, and strong Morita
    equivalence and exemplify their properties for star product
    algebras. The complete classification of Morita equivalent star
    products on general Poisson manifolds is discussed as well as the
    complete classification of covariantly Morita equivalent star
    products on a symplectic manifold with respect to some Lie algebra
    action preserving a connection.
\end{abstract} 

%
% Introduction
%

\section{Introduction}
\label{sec:Intro}

Morita theory is a classical topic on ring theory having a more
analytic cousin called strong Morita theory in the context of
$C^*$-algebras. The underlying idea is that one wants to learn
something about the categories of modules over a given ring by
comparing it with the corresponding category for some other ring: even
if the category of left modules might be a very complicated object, it
can still be possible to state that for two (non-isomorphic) rings the
corresponding categories are equivalent. In the realm of unital rings
this is precisely Morita equivalence. In the $C^*$-algebraic framework
one is interested not just in modules but in $^*$-representations on
Hilbert spaces or, more generally, on Hilbert modules over some
auxiliary $C^*$-algebra.

In deformation quantization one is interested in the representation
theories of the deformed algebras. But now the star product algebras
are more specific than ``just a ring'' and hence a purely
ring-theoretic treatment would not seem to be appropriate. It simply
will not capture all interesting properties of the star products. One
can achieve Hermitian star products yielding $^*$-algebras and using
the ring ordering of $\mathbb{R}[[\hbar]]$ one has a natural notion of
positivity at hand. Thus from this and many other aspects the star
product algebras behave much more like $C^*$-algebras. Hence one
requires a more refined notion of representation theory leading to
$^*$-representations on pre Hilbert spaces over $\mathbb{C}[[\hbar]]$
as well as on pre Hilbert modules over auxiliary $^*$-algebras over
$\mathbb{C}[[\hbar]]$. It turns out that one can transfer the notions
of strong Morita theory from the $C^*$-algebraic theory into this
entirely algebraic framework, thereby extending the previous notions
tremendously.  This way, star product algebras and $C^*$-algebras can
be treated almost on the same footing even though the star product
algebras are not at all $C^*$-algebras.

In this review, we will focus on the recent developments in the
understanding of the Morita theory of deformed algebras in general and
of the star product algebras as major class of examples. Since
recently in \cite{bursztyn.dolgushev.waldmann:2009a:pre} the final
classification of Morita equivalent star products in the general case
of Poisson manifolds was obtained, it seems to be a good point to give
such an overview.

In Section~\ref{sec:MoritaEquivalence}, we start with some elementary
presentation of the ring-theoretic aspects of Morita theory as it can
be found in any algebra textbook, see e.g.  \cite{lam:1999a}. Then we
pass to the notion of $^*$-Morita equivalence as it was developed by
Ara \cite{ara:1999a} and to the notion of strong Morita equivalence
which was first established by Rieffel for $C^*$-algebras in
\cite{rieffel:1974b} and then for general $^*$-algebras over ordered
rings in \cite{bursztyn.waldmann:2005b, bursztyn.waldmann:2001a}. Here
we stress the functorial aspects of Morita theory and discuss in
particular the functoriality of the classical limit.
Section~\ref{sec:DeformationQuantization} contains a brief
introduction to the existence and classification results in
deformation theory based on Kontsevich's formality
theorem\cite{kontsevich:2003a}. Our focus is on the notions of
equivalence of star products and of formal Poisson structures.  In
Section~\ref{sec:MoritaEquivalenceStarProducts} we explain the main
result of \cite{bursztyn.dolgushev.waldmann:2009a:pre} by first
establishing the gauge action of formal series of closed two-forms on
formal Poisson structures. On the level of equivalence classes this
provides precisely the description of Morita equivalent star products
in terms of classical data \emph{provided} the two-form is integral.
The particular case of symplectic star products is easier and was
discussed earlier in \cite{bursztyn.waldmann:2002a}. In
Section~\ref{sec:IncorporatingSymmetries} we recall the basic notions
of Morita theory in presence of symmetries which will be modelled by a
Hopf algebra action. Here the equivalence bimodules are required to
carry an action of the Hopf algebra, too, such that all structure maps
have nice covariance properties. One can now study the relations
between the various notions of Morita equivalence. Finally, in
Section~\ref{sec:gActionsSymplecticManifolds} we outline the main
results of \cite{jansen.neumaier.waldmann:2005a:pre} where the
classification of invariant star products on a symplectic manifold up
to covariant Morita equivalence with respect to a Lie algebra action
was obtained.

\noindent
\textbf{Acknowledgement:} It is a great pleasure to thank the
organisers of the Poisson 2010 conference in Rio de Janeiro, and in
particular Henrique Bursztyn, for the wonderful conference. Moreover,
I would like to thank the participants for various suggestions and
comments on my talk.

%
% Morita Equivalence
%

\section{Morita Equivalence}
\label{sec:MoritaEquivalence}

In this section we first recall some basic notions of ring-theoretic
Morita equivalence and specialize this to $^*$-algebras over ordered
rings to establish the notions of $^*$-Morita equivalence and strong
Morita equivalence.

We consider unital algebras over some fixed commutative unital ring
$\ring{C}$. Later on, in deformation quantization the cases $\ring{C}
= \mathbb{C}$ or $\mathbb{C}[[\hbar]]$ will be used. One can abandon
the condition of having unital algebras by imposing some slightly
weaker requirements (non-degeneracy and idempotency) but we do not
need this more general framework here.

The idea of Morita theory is to replace ordinary algebra homomorphisms
$\Phi\colon \mathcal{A} \longrightarrow \mathcal{B}$ by something more
general in order to have more freedom when comparing algebras. The new
arrows between algebras will now be \emph{bimodules}. Consider a
$(\mathcal{B}, \mathcal{A})$-bimodule $\BEA$: this notation says that
$\mathcal{B}$ acts from the left while $\mathcal{A}$ acts from the
right. By convention, all bimodules will have a compatible
$\ring{C}$-module structure and all structure maps will be
$\ring{C}$-(multi-)linear in the following. Moreover, the units of
$\mathcal{A}$ and $\mathcal{B}$ will act as identity on $\BEA$. We
will view $\BEA$ now as an arrow from $\mathcal{A}$ to
$\mathcal{B}$. If $\CFB$ is another bimodule for some additional
algebra $\mathcal{C}$ then the composition of bimodules will be the
tensor product over the algebra in the middle, i.e.
\begin{equation}
    \label{eq:BimoduleComposition}
    \bfig
    \morphism(0,0)|b|/{@{<-}@/^2em/}/<800,0>[\mathcal{C}`\mathcal{B};\CFB]
    \morphism(800,0)|b|/{@{<-}@/^2em/}/<800,0>[\mathcal{B}`\mathcal{A}.;\BEA]
    \morphism(0,0)|b|/{@{<-}@/_2em/}/<1600,0>[\mathcal{C}`\mathcal{A}.;\CFB
    \otimes_{\scriptscriptstyle \mathcal{B}}\BEA]
    \efig
\end{equation}
As a unit morphism one uses the canonical bimodule $\AAA$ where
$\mathcal{A}$ acts on itself by left and right multiplications. Now
$\tensor$ is not directly associative but only associative up to a
natural isomorphism. Also, $\AAA$ is not directly a unit element for
$\tensor$ but again only up to a natural isomorphism. Thus we have to
identify isomorphism classes of bimodules: this results in an honest
category which we denote by $\BiMod$. Strictly speaking, we should
confine ourselves to bimodules and algebras from some Grothendieck
universe in order to get honest sets of morphisms. But this is a
technical issue which will not affect the notion of Morita equivalence
at all.

While passing to isomorphism classes of bimodules has the advantage to
yield a category we can also stay with the bimodules directly and
taking into account that $\tensor$ is not really associative. This
leads to a \emph{bicategory} where the $1$-morphisms are the
bimodules and the $2$-morphisms, i.e. the morphisms between the
$1$-morphisms, are the bimodule homomorphisms. Such $2$-morphisms can
then be depicted by
\begin{equation}
    \label{eq:TwoMorphism}
    \bfig
    \morphism(0,0)|a|/{@{<-}@/^2em/}/<800,0>[\mathcal{B}`\mathcal{A};\BEA]
    \morphism(0,0)|b|/{@{<-}@/_2em/}/<800,0>[\mathcal{B}`\mathcal{A};\BEpA]
    \morphism(400,150)|b|/=>/<0,-300>[`;\Phi]
    \efig
\end{equation}
for a bimodule homomorphism $\Phi: \BEA \longrightarrow \BEpA$. Note
that we require $\Phi$ to be $\ring{C}$-linear as well. It is then a
classical result that the natural isomorphisms implementing the
associativity for $\tensor$ and the unit properties of $\AAA$ satisfy
the necessary coherence properties to yield a bicategory (weak
$2$-category), see \cite{benabou:1967a}. This bicategory will then be
denoted by $\BiModBI$.

We are now in the position to state the definition of Morita
equivalence:
\begin{definition}[Morita equivalence and Picard group]
    \label{definition:MoritaEquivalence}%
    Two unital $\ring{C}$-algebras $\mathcal{A}$ and $\mathcal{B}$ are
    Morita equivalent if they are isomorphic in $\BiMod$. A bimodule
    $\BEA$ representing an invertible arrow in $\BiMod$ is called an
    equivalence bimodule. The groupoid of invertible arrows in
    $\BiMod$ is called the Picard groupoid, denoted by $\Pic$. The
    isotropy group of it at $\mathcal{A}$ is called the Picard group
    $\Pic(\mathcal{A})$ of $\mathcal{A}$.
\end{definition}
Alternatively, we can also use isomorphisms in $\BiModBI$ in the sense
of bicategories, i.e. two objects $\mathcal{A}$ and $\mathcal{B}$ are
isomorphic if there are $1$-morphisms in both directions such that
their compositions are isomorphic (via $2$-morphisms) to the identity
morphisms $\AAA$ and $\BBB$, respectively. This gives then the
\emph{Picard bigroupoid} $\PicBI$ as well as the \emph{Picard bigroup}
$\PicBI(\mathcal{A})$ at $\mathcal{A}$.

The main task of Morita theory is then twofold: first one would like
to know which algebras are Morita equivalent, this is described by the
\emph{orbits} of the Picard groupoid. Second, one would like to
understand in how many different ways two algebras can be Morita
equivalent. Thanks to the groupoid structure this is equivalent to
determine the isotropy groups, i.e. the Picard groups: they encode how
many self-equivalences an algebra has. The classical theorem of Morita
determines now the structure of equivalence bimodules:
\begin{theorem}[Morita]
    \label{theorem:Morita}%
    $\BEA$ is an equivalence bimodule for two unital
    $\ring{C}$-algebras $\mathcal{A}$ and $\mathcal{B}$ iff there is
    an idempotent $e = e^2 \in M_n(\mathcal{A})$ in some matrix
    algebra over $\mathcal{A}$ such that the right
    $\mathcal{A}$-module $\EA$ is isomorphic to $e \mathcal{A}^n$ and
    one has $\mathcal{A}e\mathcal{A} = \mathcal{A}$ as well as
    $\mathcal{B} \cong \End_{\mathcal{A}}(\EA) \cong e
    M_n(\mathcal{A}) e$.
\end{theorem}
Here $\mathcal{A}e\mathcal{A}$ denotes the two-sided ideal generated
by the $n^2$ components of $e$ and the isomorphism in the last
statement are those induced by the action of $\mathcal{B}$ and
$eM_n(\mathcal{A})e$ on $\EA \cong e\mathcal{A}^n$. In particular, an
equivalence bimodule $\BEA$ is a finitely generated and projective
module over $\mathcal{A}$. This shows already that we have to look for
equivalence bimodules inside the $K_0$-theory of $\mathcal{A}$.  Since
Morita equivalence is a symmetric relation by the very definition, we
can equivalently formulate things in terms of $\mathcal{B}$ instead.
\begin{example}
    \label{example:SerreSwan}%
    As a first example relevant for the following we consider the
    smooth complex-valued functions $\mathcal{A} = C^\infty(M)$ for a
    smooth manifold $M$. Then it is well-known that any finitely
    generated projective module over $C^\infty(M)$ is isomorphic as a
    right module to the sections $\Gamma^\infty(E)$ of a complex
    vector bundle $E \longrightarrow M$. Moreover, this gives an
    equivalence bimodule with $\End_{C^\infty(M)}(\Gamma^\infty(E))
    \cong \Gamma^\infty(\End(E))$ iff the fiber dimension of $E$ is
    not zero. Thus the Morita equivalent algebras to $C^\infty(M)$ are
    isomorphic to the sections $\Gamma^\infty(\End(E))$ of the
    endomorphism bundle $\End(E)$ for arbitrary non-zero vector
    bundles $E \longrightarrow M$.
\end{example}

This example also allows to determine the Picard group of
$C^\infty(M)$: Since the only way to get the endomorphism algebra
$\Gamma^\infty(\End(E))$ to be isomorphic to $C^\infty(M)$ is by a
\emph{line bundle} $L \longrightarrow M$ we have the following result:
Implementing the isomorphism by choosing an appropriate automorphism
of $C^\infty(M)$ we can arrange to get a \emph{symmetric} bimodule
where $C^\infty(M)$ acts on $\Gamma^\infty(L)$ from left and right in
the same way. Since $\Gamma^\infty(L)$ determines $L$ completely, we
are left with the classification of line bundles, which is done via
the Chern class. Since the automorphisms of $C^\infty(M)$ are just the
pull-backs with diffeomorphism, we arrive at the result that
\begin{equation}
    \label{eq:PicardCinftyM}
    \Pic(C^\infty(M))
    = \mathrm{Diffeo}(M) \ltimes \mathrm{H}^2(M, \mathbb{Z})
\end{equation}
as group where the semidirect product comes from the usual action of
the diffeomorphism on the integral cohomology classes by pull-backs.

Let us now pass to the more specific case where the underlying scalars
are of the form $\ring{C} = \ring{R}(\I)$ with an ordered ring
$\ring{R}$ and $\I^2 = -1$. In this case we have on one hand the
complex conjugation in $\ring{C}$ and on the other hand, inherited
from the ordering of $\ring{R}$, the notion of positivity. We want to
transfer this now to algebras over $\ring{C}$ as well: instead of
general unital algebras we consider now $^*$-algebras, i.e. algebras
equipped with an anti-linear, involutive anti-automorphism, the
\emph{$^*$-involution} denoted by $a \mapsto a^*$. Then we can speak
of positivity in the following way: a linear functional $\mathcal{A}
\longrightarrow \ring{C}$ is \emph{positive} if $\omega(a^*a) \ge 0$
for all $a \in \mathcal{A}$. Using this, we say that an algebra
element $a \in \mathcal{A}$ is \emph{positive} if $\omega(a) \ge 0$
for all positive linear functionals. We denote the positive algebra
elements by $\mathcal{A}^+$. Clearly, $a^*a$ and any convex
combination of such elements are in $\mathcal{A}^+$ but there might be
more. Note also that there are other scenarios where one employs a
more sophisticated version of positivity for the price of additional
structures like for $O^*$-algebras, see e.g.\cite{schmuedgen:1990a}.

As before we want now to replace the obvious notion of
$^*$-homomorphism by some bimodule version. Here we rely on the
particular case of $C^*$-algebras where this theory was studied
first. However, the essence is entirely algebraic and thus works in
our general setting as well. We consider again a $(\mathcal{B},
\mathcal{A})$-bimodule $\BEA$, now together with a \emph{inner
  product}
\begin{equation}
    \label{eq:InnerProduct}
    \SPA{\cdot, \cdot}: \BEA \times \BEA \longrightarrow \mathcal{A},
\end{equation}
such that $\SPA{\cdot, \cdot}$ is $\ring{C}$-linear in the second
argument and satisfies $\SPA{x, y \cdot a} = \SPA{x, y} a$ as well as
$\SPA{x, y} = (\SPA{y, x})^*$, and $\SPA{b \cdot x, y} = \SPA{x, b^*
  \cdot y}$ for all $x, y \in \BEA$, $a \in \mathcal{A}$, and $b \in
\mathcal{B}$. Finally, we require $\SPA{\cdot, \cdot}$ to be
non-degenerate. In this case we will call $\BEA$ together with
$\SPA{\cdot, \cdot}$ an \emph{inner product $(\mathcal{B},
  \mathcal{A})$-bimodule}. Note that the definition is not symmetric
in $\mathcal{A}$ and $\mathcal{B}$. As a second variant we consider a
\emph{completely positive} inner product where we require in addition
\begin{equation}
    \label{eq:CP}
    \left(\SPA{x_i, x_j}\right) \in M_n(\mathcal{A})^+.
\end{equation}
Here $n \in \mathbb{N}$ and $x_1, \ldots, x_n \in \BEA$ and the
$^*$-algebra structure on $M_n(\mathcal{A})$ is the usual one induced
from the one on $\mathcal{A}$. In general it will be quite difficult
to determine whether an inner product is completely positive because
we first have to determine all positive linear functionals of
$M_n(\mathcal{A})$ in order to determine the positivity of a
matrix. However, there are many examples and more particular
situations where things simplify. A bimodule with such a completely
positive inner product is then called a \emph{pre Hilbert bimodule}.

The inner product bimodules as well as pre Hilbert bimodules are now
used to define (bi-) categories $\starBiMod$ and $\StrBiMod$
($\starBiModBI$ and $\StrBiModBI$, respectively) by adapting the
notion of the tensor product appropriately. In fact, it is rather easy
to see that the definition
\begin{equation}
    \label{eq:RieffelFormula}
    \SPFEA{\phi \tensor x, \psi \tensor y}
    =
    \SPEA{x, \SPFB{\phi, \psi} \cdot y}
\end{equation}
for $x, y \in \BEA$ and $\phi, \psi \in \CFB$ extends to a
well-defined inner product on $\CFB \tensor[\mathcal{B}] \BEA$
\emph{except} that it might still be degenerate. However, a further
quotient by the degeneracy space will yield a well-defined inner
product. The quotient together with this new inner product then
defines a new composition, the tensor product $\itensor[\mathcal{B}]$.
It is slightly more tricky to see that also complete positivity is
preserved by $\itensor[\mathcal{B}]$. In both cases, with this new
tensor product we get honest categories $\starBiMod$ and $\StrBiMod$
after passing to isometric isomorphism classes or, without
identifying, bicategories $\starBiModBI$ and $\StrBiModBI$,
respectively. In the bicategory case we have to require that the
$2$-morphisms are not just bimodule morphisms but \emph{adjointable}
bimodule morphisms, i.e. there exists an adjoint with respect to the
inner product. In analogy to
Definition~\ref{definition:MoritaEquivalence} we can now state the
following:
\begin{definition}[$^*$-Morita and strong Morita equivalence]
    \label{definition:StrongStarME}%
    Two unital $^*$-algebras $\mathcal{A}$ and $\mathcal{B}$ are
    $^*$-Morita equivalent or strongly Morita equivalent if they are
    isomorphic in $\starBiMod$ or in $\StrBiMod$, respectively. A
    bimodule $\BEA$ representing an invertible arrow in $\starBiMod$
    or $\StrBiMod$ is called a $^*$-equivalence or strong equivalence
    bimodule, respectively. The groupoid of invertible arrows in
    $\starBiMod$ and in $\StrBiMod$ is called the $^*$-Picard groupoid
    $\starPic$ and the strong Picard groupoid $\StrPic$. The isotropy
    groups of them at $\mathcal{A}$ are called the $^*$-Picard group
    $\starPic(\mathcal{A})$ and the strong Picard group
    $\StrPic(\mathcal{A})$ of $\mathcal{A}$.
\end{definition}

For unital $^*$-algebras one can prove that the tensor product
$\itensor$ of $^*$-equivalence or strong equivalence bimodules does
not require the additional quotient procedure. Moreover, it is easy to
see that after forgetting about the inner product one obtains an
equivalence bimodule in the ring-theoretic sense. This yields
well-defined groupoid morphisms
\begin{equation}
    \label{eq:PicMorphisms}
    \bfig
    \Vtriangle(0,0)/>`>`>/<600,300>[\StrPic`\starPic`\Pic;``]
    \efig
\end{equation}
by successively forgetting structure. On the level of the Picard
groups many properties of these forgetful morphisms have been
discussed in \cite{bursztyn.waldmann:2005b}. In particular, even for a
$C^*$-algebra $\mathcal{A}$ over $\mathbb{C}$ the group morphism
$\StrPic(\mathcal{A}) \longrightarrow \Pic(\mathcal{A})$ is in general
nor surjective, though always injective.

In a last step we consider now formal deformations. Here we will
consider a unital algebra $\mathcal{A}$ over $\ring{C}$ together with
a formal associative deformation $\star$ in the sense of Gerstenhaber
\cite{gerstenhaber:1964a}. This means we have $\ring{C}$-bilinear maps
$C_r\colon \mathcal{A} \times \mathcal{A} \longrightarrow \mathcal{A}$
such that
\begin{equation}
    \label{eq:FormalDeformation}
    a \star b = \sum_{r=0}^\infty \hbar^r C_r(a, b)
\end{equation}
defines a $\ring{C}[[\hbar]]$-bilinear associative product for
$\mathcal{A}[[\hbar]]$ where $C_0(a, b) = ab$ is the original product
of $\mathcal{A}$. Since we work in the unital setting, we require that
$\Unit_{\mathcal{A}}$ is still a unit with respect to $\star$.  We
will abbreviate this by $\deform{\mathcal{A}} = (\mathcal{A}[[\hbar]],
\star)$.

Passing to formal power series gives rise to various classical limit
maps, where we set $\hbar = 0$. On the level of algebra elements we
write
\begin{equation}
    \label{eq:ClassicalLimit}
    \cl(a) = a_0
    \quad
    \textrm{for}
    \quad
    a = \sum_{r=0}^\infty \hbar^r a_r \in \mathcal{A}[[\hbar]],
\end{equation}
which is a $\ring{C}$-linear algebra morphism $\cl\colon
(\mathcal{A}[[\hbar]], \star) \longrightarrow \mathcal{A}$.  If we
have two deformations $\deform{\mathcal{B}}$ and
$\deform{\mathcal{A}}$ of $\mathcal{B}$ and $\mathcal{A}$,
respectively, then for a $(\deform{\mathcal{B}},
\deform{\mathcal{A}})$-bimodule $\deform{\BEA}$ we define the
classical limit as the quotient
\begin{equation}
    \label{eq:ClassicalLimitBimodule}
    \cl\colon \deform{\BEA}
    \longrightarrow
    \cl(\deform{\BEA}) = \deform{\BEA} \big/ \hbar \deform{\BEA}.
\end{equation}
The result $\cl(\deform{\BEA})$ is viewed as a $(\mathcal{B},
\mathcal{A})$-bimodule. It is now easy to see that $\cl(\deform{\AAA})
\cong \AAA$ and that the tensor product is compatible with $\cl$ up to
natural isomorphisms. This simple observation can be summarized as
follows. To keep track of the ring we denote the category of bimodule
over algebras over $\ring{C}$ and $\ring{C}[[\hbar]]$ by
$\BiMod_{\ring{C}}$ and $\BiMod_{\ring{C}[[\hbar]]}$,
respectively. Then one has the sub-category $\deform{\BiMod} \subseteq
\BiMod_{\ring{C}[[\hbar]]}$ of those algebras over $\ring{C}[[\hbar]]$
which are formal deformations of algebras over $\ring{C}$ as
above. For the morphisms in $\deform{\BiMod}$ we allow \emph{all}
bimodules and not just those of the form $\deform{\mathcal{E}} =
\mathcal{E}[[\hbar]]$. Then $\cl$ induces a functor
\begin{equation}
    \label{eq:ClassicalLimitFunctor}
    \cl\colon
    \deform{\BiMod} \longrightarrow \BiMod_{\ring{C}},
\end{equation}
called the classical limit functor. Hence we also get immediately a
groupoid morphism
\begin{equation}
    \label{eq:clPic}
    \cl\colon
    \deform{\Pic} \longrightarrow \Pic_{\ring{C}},
\end{equation}
which ultimately results in a group morphism
\begin{equation}
    \label{eq:clPicardGroup}
    \cl\colon
    \Pic(\deform{\mathcal{A}}) \longrightarrow \Pic(\mathcal{A})
\end{equation}
for every deformation $\deform{\mathcal{A}}$ of $\mathcal{A}$. Again,
one is interested in understanding the properties of this classical
limit morphism \eqref{eq:clPic} and in particular the behaviour of the
Picard groups under deformation \eqref{eq:clPicardGroup}. Many results
on this have been obtained in \cite{bursztyn.waldmann:2004a}.

Also for $^*$-algebras one can define a classical limit functor
similar to \eqref{eq:ClassicalLimitFunctor}: here one considers
\emph{Hermitian} deformations which are formal deformations $\star$
such that the original $^*$-involution is still a $^*$-involution also
with respect to the deformed product $\star$. Then the quotient
procedure for the bimodules has to be modified as the resulting inner
product on the naive quotient $\deform{\BEA} \big/ \hbar
\deform{\BEA}$ will in general be degenerate. Thus we have to divide
by the degeneracy space and get a functor
\begin{equation}
    \label{eq:clStarBiMod}
    \cl\colon
    \deform{\starBiMod} \longrightarrow
    \starBiMod\nolimits_{\ring{C}},
\end{equation}
which also gives groupoid and group morphisms for the $^*$-Picard
groupoid and $^*$-Picard groups, respectively. Finally, for the strong
version one has to take care once more: the complete positivity of the
inner product on the classical limit may fail. The way out is to allow
only those Hermitian deformations which are \emph{completely positive
  deformations}, see e.g. \cite{bursztyn.waldmann:2005a} for examples.

%
% Deformation Quantization
%

\section{Deformation Quantization}
\label{sec:DeformationQuantization}

After the algebraic preliminaries on Morita theory we pass now to the
geometric situation of deformation quantization of Poisson manifolds.

Let $(M, \pi_1)$ be a Poisson manifold with Poisson tensor $\pi_1 \in
\Gamma^\infty(\Anti^2 TM)$. Then a \emph{formal deformation} of
$\pi_1$ is a formal series
\begin{equation}
    \label{eq:piSeries}
    \pi = \hbar \pi_1 + \hbar^2 \pi_2 + \cdots
    \in \hbar \Gamma^\infty (\Anti^2 TM)[[\hbar]],
\end{equation}
such that we still have the Jacobi identity $\Schouten{\pi, \pi} =
0$. Such a formal series is also called a formal Poisson tensor. The
set of all formal Poisson tensors is denoted by $\FormalPoiss(M)$ and
those with fixed first order $\pi_1$ are denoted by $\FormalPoiss(M,
\pi_1)$. A \emph{formal vector field} is a formal series $X = \hbar
X_1 + \hbar^2 X_2 + \cdots \in \hbar \Gamma^\infty(TM)[[\hbar]]$. For
a formal vector field $X$ the exponential series $\exp(\Lie_X)$ of the
Lie derivative is a well-defined operator on formal series with
coefficients in some type of tensor fields on $M$. We call
$\exp(\Lie_X)$ the \emph{formal diffeomorphism} induced by $X$.  The
Baker-Campbell-Hausdorff theorem shows that the composition of two
formal diffeomorphisms $\exp(\Lie_X)$ and $\exp(\Lie_Y)$ is again a
formal diffeomorphism $\exp(\Lie_{\mathrm{BCH}(X, Y)})$ since both $X$
and $Y$ start in first order of $\hbar$ making the BCH series
$\mathrm{BCH}(X, Y)$ a well-defined formal vector field again. This
way, we obtain the \emph{group of formal diffeomorphisms}
$\FDiffeo(M)$ acting on various tensor fields.  Indeed, one can think
of formal vector fields, formal Poisson structures, formal
diffeomorphisms, etc. as the $\infty$-jet around $\hbar = 0$ of vector
fields, Poisson structures, diffeomorphisms, etc., depending smoothly
on the parameter $\hbar$. However, we will not need this point of view
here.

Since the Schouten bracket $\Schouten{\cdot, \cdot}$ is natural with
respect to the Lie derivative, it is clear that a formal
diffeomorphism $\exp(\Lie_X)$ maps a formal Poisson structure $\pi$ to
a formal Poisson structure $\exp(\Lie_X)(\pi)$ again. Moreover, the
first order term $\pi_1$ is preserved by this action. This motives the
definition that two formal Poisson structures $\pi$ and $\pi'$ are
called \emph{equivalent} if they are in the same $\FDiffeo(M)$-orbit,
i.e. if there is a formal vector field such that
\begin{equation}
    \label{eq:EquivalentPoisson}
    \E^{\Lie_X} (\pi) = \pi'.
\end{equation}
In this case we necessarily have $\pi_1 = \pi_1'$ and we write $\pi
\sim \pi'$. The equivalence classes of this equivalence relation are
then denoted by
\begin{equation}
    \label{eq:FPoiss}
    \FPoiss(M) = \FormalPoiss(M) \big/ \FDiffeo(M)
    \quad
    \textrm{and}
    \quad
    \FPoiss(M, \pi_1) = \FormalPoiss(M, \pi_1) \big/ \FDiffeo(M).
\end{equation}
They are the moduli space for the inequivalent formal deformations of
a given Poisson structure $\pi_1$.  In general, it will be very
complicated to determine the set $\FPoiss(M, \pi_1)$ for a given
Poisson structure $\pi_1$. By abstract deformation theory one can say
that $[\pi_2]$ is a well-defined class in the second Poisson
cohomology of $\pi_1$ but it is not clear which such infinitesimal
deformations can actually be lifted to formal deformations of all
order in $\hbar$. If however, $\pi_1$ is symplectic and comes from a
symplectic form $\omega_1$ then the moduli space $\FPoiss(M, \pi_1)$ is
easily be described by the inequivalent formal deformations of
$\omega$. Here the result is
\begin{equation}
    \label{eq:SymplecticCase}
    \FPoiss(M, \pi_1) = [\omega] + \hbar \HdR^2(M, \mathbb{C})[[\hbar]],
\end{equation}
where we have (artificially) put an affine space modeled on $\HdR^2(M,
\mathbb{C})[[\hbar]]$ instead of $\HdR^2(M, \mathbb{C})[[\hbar]]$
itself, just to keep track of the symplectic form we started with.
The proofs of these facts are well-known and can e.g. be found in the
textbook \cite[Sect.~4.2.4]{waldmann:2007a}.

Let us now recall the basic notions from deformation quantization
\cite{bayen.et.al:1978a}, see also the textbook \cite{waldmann:2007a}
for a detailed exposition. On a Poisson manifold $(M, \pi_1)$ a
\emph{star product} $\star$ is a formal associative deformation of
$C^\infty(M)$ as in \eqref{eq:FormalDeformation} with the additional
requirements that the first order commutator
\begin{equation}
    \label{eq:CeinsPoisson}
    C_1(f, g) - C_1(g, f) = \I \{f, g\}
\end{equation}
gives the Poisson bracket coming from $\pi_1$. Moreover, one requires
that the $C_r$ are bidifferential operators. The set of all star
products on $M$ is sometimes denoted by $\Deformation(M)$ and the star
products quantizing the Poisson structure $\pi_1$ are then
$\Deformation(M, \pi_1)$.

If $D = \hbar D_1 + \hbar^2 D_2 + \cdots \in \hbar
\Diffop(M)[[\hbar]]$ is a formal series of differential operators we
construct analogously to \eqref{eq:EquivalentPoisson} an action on
star products: the exponential series $T = \exp(D)$ is a well-defined
formal series of differential operators, now starting with the
identity, i.e. $T = \id + \hbar T_1 + \cdots$. Conversely, every such
formal series is of this form as we can always build $D = \log(T)$ as
a well-defined formal power series. Note that $D$ vanishes on
constants iff $T$ is the identity on constants. Now if $\star$ is a
star product for $\pi_1$ then also
\begin{equation}
    \label{eq:EquivalentStarProduct}
    f \star' g = T^{-1}(Tf \star Tg)
\end{equation}
is easily shown to be a star product quantizing the same Poisson
structure $\pi_1$. Here we need that $D = \log(T)$ vanishes on
constants in order to have again $1 \star' f = f = f \star' 1$. This
allows to interpret the operators $D$ as quantum analogs of formal
vector fields while the operators $T$ are the quantum analogs of
formal diffeomorphisms. We call such an operator $T$ an
\emph{equivalence transformation}. Clearly, we get a group structure
by multiplying equivalence transformations which corresponds to the
Lie algebra structure of the operators $D$ coming from the
commutator. We end up with an action of the group of equivalence
transformations on $\Deformation(M)$ which preserves each
$\Deformation(M, \pi_1)$. This allows to define two star products
$\star$ and $\star'$ to be \emph{equivalent}, denoted by $\star \sim
\star'$, if they are in the same orbit under the action
\eqref{eq:EquivalentStarProduct} of the equivalence
transformations. This gives us the analog of \eqref{eq:FPoiss} and we
set
\begin{equation}
    \label{eq:DefM}
    \Def(M) = \Deformation(M) \big/ \sim
    \quad
    \textrm{and}
    \quad
    \Def(M, \pi_1) = \Deformation(M, \pi_1) \big/ \sim.
\end{equation}

One of the major achievements in deformation quantization is now the
famous statement of Kontsevich that the two moduli spaces $\FPoiss(M,
\pi_1)$ and $\Def(M, \pi_1)$ are in bijection for every Poisson
structure $\pi_1$. More precisely, the formality map $\mathcal{K}$ of
Kontsevich gives a construction where a formal Poisson structure $\pi
= \hbar \pi_1 + \hbar^2 \pi_2 + \cdots$ is used to build a formal star
product $\star_\pi$ in such a way that $\star_\pi \sim \star_{\pi'}$
iff $\pi \sim \pi'$. The precise construction
\begin{equation}
    \label{eq:pitostarpi}
    \pi \; \mapsto \; \star_\pi
\end{equation}
requires the formality map $\mathcal{K}$ and is involved, both from
the conceptual point of view as well as technically, see
\cite{kontsevich:2003a} for further details. Thus the choice of a
formality map results in a bijection
\begin{equation}
    \label{eq:FPoissToDef}
    \mathcal{K}_*\colon
    \FPoiss(M, \pi_1) \longrightarrow \Def(M, \pi_1).
\end{equation}

As a remark we note that for a real formal Poisson structure $\pi =
\cc{\pi}$ Kontsevich's formality on $\mathbb{R}^n$ produces a
Hermitian star product, i.e. $\cc{f \star_\pi g} = \cc{g} \star_\pi
\cc{f}$. Also the global formality map of Dolgushev has this property
\cite{dolgushev:2005a}. Finally, one can show that a Hermitian star
product is always a completely positive deformation
\cite{bursztyn.waldmann:2005a}.

%
% Morita Equivalence of Star Products
%

\section{Morita Equivalence of Star Products}
\label{sec:MoritaEquivalenceStarProducts}

The main question concerning Morita theory in deformation quantization
is now which star products are Morita equivalent. Since we have a good
understanding of the equivalence classes of star products in terms of
the equivalence classes of formal Poisson tensors one can refine the
task as follows: describe the Morita equivalence of $\star_\pi$ and
$\star_{\pi'}$ in terms of the equivalence classes of $\pi$ and
$\pi'$. Indeed, since isomorphic algebras are Morita equivalence the
Morita equivalence of $\star_\pi$ and $\star_{\pi'}$ will only depend
on the \emph{classes} of $\pi$ and $\pi'$. Moreover, since every star
product is equivalent (and hence isomorphic) to a $\star_\pi$, it
suffices to consider those.

There are now two reasons for star products to be Morita equivalent
which we would like to discuss separately. First it is clear that any
Poisson diffeomorphism $\Phi\colon (M, \pi_1) \longrightarrow (M,
\pi_1)$ maps $\star_\pi$ to an isomorphic star product
\begin{equation}
    \label{eq:PhipullbackStar}
    \star_\pi \; \mapsto \; \Phi^*(\star_\pi) \sim \star_{\Phi^*\pi},
\end{equation}
where we use in the second step that the global formalities have a
good covariance property up to equivalence. Since $\star_\pi$ and
$\Phi^*(\star_\pi)$ are isomorphic (via $\Phi^*$) they are Morita
equivalent in a trivial way. This is the simple part of the
description. Note that for star products quantizing diffeomorphic but
different Poisson structures we still can have isomorphism via general
diffeomorphisms. However, we study Morita equivalence of star products
for a fixed Poisson structure $\pi_1$ on $M$.

The non-trivial part of Morita equivalence comes from the non-trivial
classical equivalence bimodules, the line bundles. We know from the
classical limit morphism \eqref{eq:clPic} that the classical limit of
an equivalence bimodule has to be a classical equivalence
bimodule. Conversely, given a line bundle $L \longrightarrow M$ one
can show that there is a unique way up to equivalence to deform the
right module structure of $\Gamma^\infty(L)$ into a right module
structure $\bullet$ for $\Gamma^\infty(L)[[\hbar]]$ with respect to
the star product algebra $(C^\infty(M)[[\hbar]], \star)$. Moreover, it
turns out that the module endomorphisms of this new, deformed right
module are in bijection to $\Gamma^\infty(\End(L))[[\hbar]] =
C^\infty(M)[[\hbar]]$. Hence we get an induced deformed product for
$C^\infty(M)[[\hbar]]$ which turns out to be a star product
$\star'$. Being isomorphic to the module endomorphisms this induces
also a left module structure $\bullet'$ for $\star'$ such that we get
a bimodule in the end.  Finally, $\star'$ is uniquely determined by
$L$ and $\star$ up to equivalence since $\bullet$ was unique up to
equivalence. Thus we get, on the level of equivalence classes, a
well-defined map
\begin{equation}
    \label{eq:Lacts}
    L\colon [\star] \mapsto [\star'].
\end{equation}
It is now easy to see that the deformed bimodule is still a Morita
equivalence bimodule and all Morita equivalences arise this way. These
results have been obtained very early in
\cite{bursztyn.waldmann:2000b, bursztyn:2002a}.

The remaining task is now to compute $\star'$ for a given $\star =
\star_\pi$ and determine the corresponding $\pi'$ such that $\star'
\sim \star_{\pi}$.

The main idea how this is achieved is to use local transition
functions to describe $L$. Then these transition functions allow for a
suitable quantum analog obeying a cocycle identity with respect to the
star product. This gives a local description of the deformed right
module structure and hence also a local description of
$\star'$. Moreover, locally the two star products $\star'$ and $\star$
are even equivalent and the difference between equivalence and Morita
equivalence is a global effect. Next one uses a two-form $B$
representing $2\pi \I c_1(L)$, e.g. the curvature of a connection on
$L$. Then the idea is to pass from the deformed transition functions
to local expressions involving $B$. Here comes now the following
construction into the game. Recall that closed two-forms act on
Poisson structures, at least on the formal level, as follows: for $B
\in \Gamma^\infty(\Anti^2 T^*M)[[\hbar]]$ and a formal Poisson
structure $\pi \in \hbar \Gamma^\infty(\Anti^2 TM)[[\hbar]]$ we
consider the corresponding (formal) bundle maps $B^\sharp$ and
$\pi^\sharp$. Then the \emph{gauge transformation} of $\pi$ by $B$ is
defined to be the formal bivector field $\tau_B(\pi)$ characterized by
\begin{equation}
    \label{eq:tauB}
    \tau_B(\pi)^\sharp =
    \pi^\sharp \circ \frac{1}{\id + B^\sharp \circ \pi^\sharp}.
\end{equation}
Since $\pi$ starts in first order of $\hbar$ the inverse is
well-defined indeed.

If $B = \D A$ is an exact two-form then we can build a formal vector
field out of the potential $A$ via $\pi$ and it is a straightforward
computation that the corresponding formal diffeomorphism maps
$\tau_B(\pi)$ to $\pi$. Thus in the exact case the gauged bivector
field is equivalent to $\pi$. In particular, it is again a formal
Poisson structure. Since the Jacobi identity is a local property and
since closed two-forms are locally exact, $\tau_B(\pi)$ is always a
Poisson structure if $\D B = 0$. However, in general it will no longer
be equivalent to $\pi$. Note however, that the first order of
$\tau_B(\pi)$ coincides with $\pi_1$.

Another simple computation shows that the formal closed two-forms
\emph{act} on formal Poisson structures via \eqref{eq:tauB} where we
view the closed formal two-forms as abelian group with respect to the
usual addition. Then the above results show that we get a well-defined
action on the level of deRham cohomology classes on one hand and on
equivalence classes of formal Poisson structures on the other hand,
i.e.
\begin{equation}
    \label{eq:Hzweiacts}
    \HdR^2(M, \mathbb{C})[[\hbar]]
    \circlearrowleft
    \FPoiss(M)
    \quad
    \textrm{and}
    \quad
    \HdR^2(M, \mathbb{C})[[\hbar]]
    \circlearrowleft
    \FPoiss(M, \pi_1).
\end{equation}
The following statement is now the first main result of
\cite{bursztyn.dolgushev.waldmann:2009a:pre}. The formula was already
found in \cite{jurco.schupp.wess:2002a} on more heuristic arguments:
\begin{theorem}
    \label{theorem:Classification}%
    Let $(M, \pi_1)$ be a Poisson manifold and $\pi$ a formal Poisson
    structure with first order $\pi_1$. Let $L \longrightarrow M$ be a
    line bundle with $B \in \Gamma^\infty(\Anti^2 T^*M)$ representing
    $2 \pi \I c_1(L)$. Then the star product $\star'$ obtained from
    \eqref{eq:Lacts} out of $\star_\pi$ is equivalent to
    $\star_{\pi'}$ with
    \begin{equation}
        \label{eq:piPrime}
        \pi' = \tau_B (\pi).
    \end{equation}
\end{theorem}
This way, one has the full classification of star products up to
ring-theoretic Morita equivalence: it only remains to take into
account the simpler part, i.e. the Poisson diffeomorphisms as
discussed above. It should be noted that the classification of star
products by formal deformations of the Poisson structure $\pi_1$
requires the \emph{choice} of a global formality. It should also be
noted that in the proof one needs a formality which is differential
and vanishes on constants. The global formality constructed in
\cite{dolgushev:2005a} fulfills all these requirements.
\begin{remark}
    \label{remark:SymplecticCase}%
    In the case of symplectic manifolds one has an alternative
    classification of star products by means of their
    \emph{characteristic class} $c(\star) \in \frac{[\omega]}{\I\hbar}
    + \HdR^2(M, \mathbb{C})[[\hbar]]$, see
    e.g. \cite{gutt.rawnsley:1999a}. Here $\star$ and $\star'$ are
    equivalent iff $c(\star) = c(\star')$. The important feature is
    that the class $c$ can be defined intrinsically without reference
    to any construction method for the star products. In particular,
    it does not rely on the choice of a formality. The second main
    result of \cite{bursztyn.dolgushev.waldmann:2009a:pre} is that the
    characteristic class $c(\star_\pi)$ of $\star_\pi$ constructed
    from a global formality for a formal Poisson structure deforming
    the symplectic $\pi_1$ is the ``inverse'' of the class of
    $\pi$. This matches the earlier result from
    \cite{bursztyn.waldmann:2002a} that in the symplectic case one has
    \begin{equation}
        \label{eq:SymplecticFormula}
        c(\star') = c(\star) + 2 \pi \I c_1(L)
    \end{equation}
    where again $\star'$ is the star product from \eqref{eq:Lacts}.
\end{remark}
\begin{remark}
    \label{remark:StarStrongForStars}%
    As a last remark here we note that the $^*$-Morita and the strong
    Morita theory of star products is now fairly easy: from general
    results in \cite{bursztyn.waldmann:2005b} one knows that Hermitian
    star products are Morita equivalent iff they are strongly Morita
    equivalent. Moreover, the kernel and image of the groupoid
    morphisms \eqref{eq:PicMorphisms} are very well understood for the
    case of star product algebras. Thus the additional requirements of
    having (completely positive) inner products on the equivalence
    bimodules do not cause any further difficulties.
\end{remark}

%
% Incorporating Symmetries
%

\section{Incorporating Symmetries}
\label{sec:IncorporatingSymmetries}

In many applications one is not just interested in Morita equivalence
but in Morita equivalence compatible with some additional symmetry:
Here several possibilities arise like a symmetry of a Lie algebra
$\lie{g}$ acting on all algebras by derivations or a symmetry of a
group $G$ acting by automorphisms. To combine these notions it is
convenient to consider a Hopf algebra symmetry.

Let $H$ be a Hopf algebra over $\ring{C}$ which should encode the type
of symmetry which we want to study. Then we consider $H$-module
algebras, i.e. algebras $\mathcal{A}$ over $\ring{C}$ which carry an
action of $H$: there is a left module structure of $H$ on
$\mathcal{A}$ denoted by
\begin{equation}
    \label{eq:Action}
    \acts\colon H \times \mathcal{A} \longrightarrow \mathcal{A},
\end{equation}
such that in addition one has $g \acts (ab) = (g_\sweedler{1} \acts
a)(g_\sweedler{2} \acts b)$ for all $g \in H$ and $a, b \in
\mathcal{A}$. Here we use the Sweedler notation for the coproduct
$\Delta(g) = g_\sweedler{1} \tensor g_\sweedler{2}$ as usual. In the
case of $^*$-algebras we require a Hopf $^*$-algebra and the action
should fulfill $(g \acts a)^* = S(g)^* \acts a^*$ where $S$ is the
antipode of $H$. If $H = \mathcal{U}(\lie{g})$ is the universal
enveloping algebra of a Lie algebra then $\acts$ corresponds just to a
Lie algebra action by derivations and if $H = \ring{C}[G]$ is the
group algebra of some group $G$, then $\acts$ reduces to a group
representation by automorphisms. Thus we cover the two cases mentioned
above.

In a next step we consider bimodules. On a $(\mathcal{B},
\mathcal{A})$-bimodule $\BEA$ we want to implement also an action of
$H$. We require that $\BEA$ is a left $H$-module such that in addition
\begin{equation}
    \label{eq:ActsOnBimodule}
    g \acts (b \cdot x)
    =
    (g_\sweedler{1} \acts b) \cdot
    (g_\sweedler{2} \acts x)
    \quad
    \textrm{and}
    \quad
    g \acts (x \cdot a)
    =
    (g_\sweedler{1} \acts x) \cdot
    (g_\sweedler{2} \acts a)
\end{equation}
for all $g \in H$, $b \in \mathcal{B}$, $x \in \BEA$, and $a \in
\mathcal{A}$. In this case we call the bimodule $H$-covariant (or
$H$-equivariant).

It is now a simple check that $\AAA$ with its induced left $H$-module
structure is a $H$-covariant $(\mathcal{A},
\mathcal{A})$-bimodule. Moreover, the tensor product $\tensor$ of
$H$-covariant bimodules gives again a $H$-covariant bimodule by
defining the $H$-action on the tensor product according to
\begin{equation}
    \label{eq:ActsTensor}
    g \acts (\phi \tensor x)
    =
    (g_\sweedler{1} \acts \phi)
    \tensor
    (g_\sweedler{2} \acts x)
\end{equation}
for $g \in H$, $\phi \in \CFB$, and $x \in \BEA$. Using now only
$H$-equivariant bimodule morphisms to relate $H$-covariant bimodules
we obtain a category $\BiModH$ of unital $\ring{C}$-algebras with
$H$-action as objects and isomorphism classes of $H$-covariant
bimodules as morphisms. Not yet identifying bimodules up to
isomorphisms yields again a bicategory, denoted by $\BiModBIH$.

There is also a way to incorporate inner products. If $\BEA$ is a
inner product $(\mathcal{B}, \mathcal{A})$-bimodule with an $H$-action
then the compatibility with the inner product we need is
\begin{equation}
    \label{eq:ActsSP}
    g \acts \SPA{x, y}
    =
    \SPA{S(g_\sweedler{1})^* \acts x, g_\sweedler{2} \acts y}
\end{equation}
for all $g \in H$ and $x, y \in \BEA$. In this case we call $\BEA$ a
$H$-covariant inner product bimodule. For the morphisms between
$H$-covariant inner product bimodules we take now the $H$-equivariant
adjointable bimodule morphisms. Again, it is a routine check that
$\itensor$ is compatible with this additional symmetry.  This gives
the category $\starBiModH$ with objects being the unital $^*$-algebras
over $\ring{C}$ with a $^*$-action of $H$ and isometric
$H$-equivariant isomorphism classes of $H$-covariant inner product
bimodules as morphisms. If in addition we take completely positive
inner product, no further compatibility is needed. We obtain the
category $\StrBiModH$. Again, we also have bicategories
$\starBiModBIH$ and $\StrBiModBIH$ in this case, the check of the
needed coherences is slightly more involved but still
straightforward. Details on this can be found in
\cite{jansen.waldmann:2006a, jansen:2006a} as well as in
\cite{calon:2010a}.

It is now clear how to define Morita equivalence in the presence of
symmetries: we use isomorphism in the categories $\BiModH$,
$\starBiModH$, and $\StrBiModH$ to define \emph{$H$-covariant Morita
  equivalence}, \emph{$H$-covariant $^*$-Morita equivalence} and
\emph{$H$-covariant strong Morita equivalence}, respectively.  The
corresponding groupoids of invertible arrows in these categories are
then the \emph{$H$-covariant Picard groupoid} $\PicH$, the
$H$-covariant $^*$-Picard groupoid $\starPicH$, and the $H$-covariant
strong Picard groupoid $\StrPicH$, respectively.

In \cite{jansen.waldmann:2006a, jansen.neumaier.waldmann:2005a:pre,
  calon:2010a} many general statements about these notions of
$H$-covariant Morita equivalence have been established. We mention
just three aspects relevant for deformation quantization:

First one can define again a classical limit functor. Here one can
even allow for a formal Hopf algebra deformation $\deform{H}$ of $H$
with deformed structure maps as usual. Then we get again a
sub-category $\deform{\BiModH}$ of $\BiMod_{\deform{H}}$ (defined over
$\ring{C}[[\hbar]]$ as before) with objects being those unital
algebras $\deform{\mathcal{A}}$ over $\ring{C}[[\hbar]]$ with
$\deform{H}$-action which are deformations of unital algebras
$\mathcal{A}$ over $\ring{C}$ with $H$-action. Then we get a classical
limit functor
\begin{equation}
    \label{eq:clBiModH}
    \cl\colon
    \deform{\BiModH} \longrightarrow \BiModH,
\end{equation}
restricting to a groupoid morphism on the level of the corresponding
covariant Picard groupoids, eventually leading to a group morphism
\begin{equation}
    \label{eq:clPicH}
    \cl\colon
    \Pic_{\deform{H}}(\deform{\mathcal{A}}) \longrightarrow \PicH(\mathcal{A}).
\end{equation}
Analogously, there is a $^*$-version for Hermitian deformations and a
strong version for completely positive deformations as well. We do not
spell out the detail which should be clear. First steps in
understanding the kernel of \eqref{eq:clPicH} are available in
\cite{calon:2010a}. It should also be mentioned that the classical
limit is already available on the level of the bicategories and gives
a homomorphism of bicategories there (a weak form of a $2$-functor).

Second, one can successively forget information. This gives groupoid
morphisms leading to the commuting diagram
\begin{equation}
    \label{eq:BigDiagram}
    \bfig
    \Vtriangle(0,450)<800,300>[\StrPicH`\starPicH`\Pic_H;``]
    \Vtriangle(0,0)/@{>}|\hole`>`>/<800,300>[\StrPic`\starPic,`\Pic;``]
    \morphism(0,650)<0,-250>[`;]
    \morphism(800,350)<0,-250>[`;]
    \morphism(1600,650)<0,-250>[`;]
    \efig
\end{equation}
between the various types of Picard groupoids. Additionally, the
classical limit fits into this diagram nicely and gives yet some more
groupoid morphisms making the doubled diagram also commutative. Again,
kernels and images of these forgetful morphisms have been studied in
various contexts.

Third, one can study the Morita equivalence of crossed products using
$H$-covariant Morita equivalence. Having an $H$-action on
$\mathcal{A}$ gives us a crossed product algebra structure on the
tensor product $\mathcal{A} \tensor H$ which we shall denote by
$\mathcal{A} \rtimes H$. Moreover, for a $H$-covariant bimodule $\BEA$
one can construct a $(\mathcal{B} \rtimes H, \mathcal{A} \rtimes
H)$-bimodule structure on the tensor product $\mathcal{E} \tensor H$,
the resulting bimodule will then be denoted by $\mathcal{E} \rtimes
H$. This results in a functor
\begin{equation}
    \label{eq:CrossedFunctor}
    \BiModH \longrightarrow \BiMod
\end{equation}
restricting to a groupoid morphism $\PicH \longrightarrow \Pic$ and
ultimately to a group morphism
\begin{equation}
    \label{eq:PicHCrossedPic}
    \PicH(\mathcal{A}) \longrightarrow \Pic(\mathcal{A} \rtimes H).
\end{equation}
The same construction goes through in the case of inner product
bimodules and pre Hilbert bimodules, see \cite{jansen.waldmann:2006a}
for further details and examples.

%
% $\lie{g}$-Actions on Symplectic Manifolds
%

\section{$\lie{g}$-Actions on Symplectic Manifolds}
\label{sec:gActionsSymplecticManifolds}

In this last section we apply our considerations on symmetries to the
particular case of star products on symplectic manifolds which are
invariant under a Lie algebra action of some finite-dimensional real
Lie algebra $\lie{g}$.

Thus let $(M, \omega)$ be a symplectic manifold endowed with an action
of $\lie{g}$, e.g. coming from an action of a Lie group $G$ with
$\lie{g}$ as its Lie algebra. For technical reasons we require the
action to preserve a connection $\nabla$. This is in fact a rather
mild requirement: if a Lie group $G$ acts \emph{properly} it preserves
a connection, but there are other examples of non-proper actions which
still preserve a connection as e.g. the linear action of
$\mathrm{Sp}(2n, \mathbb{R})$ on $\mathbb{R}^{2n}$ which preserves the
canonical flat connection. Without restriction we can assume that
$\nabla$ is torsion-free and symplectic in addition.

It is now a well-known fact that Fedosov's construction of a star
product, see e.g. \cite{fedosov:1996a}, gives an invariant star
product provided one has invariant geometric entrance data. More
precisely, in the situation where one has an invariant connection the
moduli space $\Def^{\lie{g}}(M, \omega)$ of invariant star products
$\Deformation^{\lie{g}}(M, \omega)$ modulo invariant equivalences is
in bijection to the second \emph{invariant} deRham cohomology
$\HdR^2(M, \mathbb{C})^{\lie{g}}[[\hbar]]$. In fact, there is a
$\lie{g}$-invariant characteristic class
\begin{equation}
    \label{eq:InvariantClass}
    c^{\lie{g}}\colon
    \Deformation^{\lie{g}}(M, \omega)
    \ni \star
    \; \mapsto \;
    c^{\lie{g}}(\star) \in
    \frac{[\omega]}{\I\hbar}
    + \HdR^2(M, \mathbb{C})^{\lie{g}}[[\hbar]],
\end{equation}
inducing this bijection. Note that under the canonical forgetful map
$\HdR^2(M, \mathbb{C})^{\lie{g}} \longrightarrow \HdR^2(M,
\mathbb{C})$ the class $c^{\lie{g}}(\star)$ is mapped to
$c(\star)$. The proofs of this result and also the case of Lie group
actions can be found in \cite{bertelson.bieliavsky.gutt:1998a}.

We are now interested in the $\lie{g}$-covariant Morita theory of
invariant symplectic star products. In more detail, we would like to
have a refined statement analogously to \eqref{eq:SymplecticFormula}
using the invariant characteristic class $c^{\lie{g}}(\cdot)$ instead.

First we study the classical limit of a $\lie{g}$-covariant
equivalence bimodule $\deform{\mathcal{L}}$ between two
$\lie{g}$-invariant star products $\star'$ acting from the left and
$\star$ acting from the right. We know that there is a unique
symplectomorphism $\Psi$ such that the twisted equivalence bimodule
${}^\Phi\deform{\mathcal{L}}$ for $\Phi^*(\star')$ and $\star$ has a
symmetric equivalence bimodule as classical limit, i.e. the smooth
sections of a line bundle $L$, unique up to isomorphism, with the same
action of $C^\infty(M)$ from the left and the right. One can show that
$\Psi$ is necessarily $\lie{g}$-equivariant and thus $\Psi^*(\star')$
is again a $\lie{g}$-invariant star product. In a next step one shows
that on the line bundle $L$ we have a lift of the $\lie{g}$-action on
$M$ together with an invariant connection $\nabla^L$: the connection
is obtained as the difference of the left and right multiplications in
the first order of $\hbar$. In fact, one can show by this construction
that the original $\deform{\mathcal{L}}$ was isomorphic to the
$\Psi$-twist of a $\lie{g}$-covariant \emph{deformation} of $L$ with
respect to two $\lie{g}$-invariant star products $\star''$ and $\star$
such that $\star''$ is $\lie{g}$-equivariantly equivalent to
$\Psi^*(\star')$. The usage of $\star''$ allows to achieve that
$\star''$ and $\star$ have the \emph{same} first order term (and not
just the same first order commutator).

In the situation of Hermitian star products and a covariant
$^*$-equivalence bimodule one has also a $\lie{g}$-invariant pseudo
Hermitian fiber metric $h_0$ on $L$ and $\nabla^L$ is metric such that
$h_0$ is the classical limit of the inner product on
$\deform{\mathcal{L}}$. Finally, in the strong equivalence case, $h_0$
is positive, i.e. a Hermitian fiber metric.

Conversely, given a line bundle $L$ over $M$ which allows for a lift
of the $\lie{g}$-action and given a $\lie{g}$-invariant connection, we
can always deform the classical bimodule structure into a quantum
bimodule structure \emph{preserving} the $\lie{g}$-invariance. This is
a simple application of Fedosov's construction adapted to vector
bundles, see \cite{waldmann:2002b}. In the Hermitian case we can also
deform a classical and $\lie{g}$-invariant fiber metric, preserving
its positivity. Since in the Fedosov construction it is very easy to
keep track of the (invariant) characteristic classes
$c^{\lie{g}}(\cdot)$, one arrives at the following theorem
\cite{jansen.neumaier.waldmann:2005a:pre}:
\begin{theorem}
    \label{theorem:EquivariantMorita}%
    Let $M$ be a symplectic manifold with $\lie{g}$-action such that a
    $\lie{g}$-invariant (symplectic torsion-free) connection of $M$
    exists. Moreover, let $\star'$ and $\star$ be two
    $\lie{g}$-invariant (Hermitian) star products on $M$. Then $\star$
    and $\star'$ are $\lie{g}$-covariantly (strongly) Morita
    equivalent iff there is a $\lie{g}$-equivariant symplectomorphism
    $\Psi$ such that $\Psi^* c^{\lie{g}}(\star') - c^{\lie{g}}(\star)$
    is in the image of the first map in 
    \begin{equation}
        \label{eq:CanonicalMapHZwei}
        \mathrm{H}^2_{\lie{g}} (M, \mathbb{C})
        \longrightarrow
        \HdR^2(M, \mathbb{C})^{\lie{g}}
        \longrightarrow
        \HdR(M, \mathbb{C})^{\lie{g}},
    \end{equation}
    and maps to a $2\pi\I$-integral class under the second map.
\end{theorem}
Here $\mathrm{H}^2_{\lie{g}}(M, \mathbb{C})$ denotes the
$\lie{g}$-equivariant deRham cohomology where we use the Cartan model
to define this cohomology. In particular, we do neither require a Lie
group action integrating the Lie algebra action of $\lie{g}$ nor any
compactness assumptions.

In fact, one even knows that the equivalence bimodules are obtained
from $\lie{g}$-equivariant symplectomorphism $\Psi$ on one hand and
from $\lie{g}$-invariant deformations of line bundles $L$ with a
lifted action on the other hand. Note also, that with the usual
arguments it is now easy to lift from the infinitesimal symmetry by
$\lie{g}$ to the integrated symmetry by the connected and
simply-connected Lie group $G$ integrating $\lie{g}$. Finally, much of
the above arguments can be carried through also in the case of a
discrete group symmetry, only the existence of invariant classical
objects is more involved in this case.

Going beyond the symplectic case to the general Poisson case should be
possible by using the equivariant formality theorem of Dolgushev
\cite{dolgushev:2005a} which are based on the existence of an
invariant connection on $M$: it is again the same classical data one
has to invest and whose existence is well-studied under various
circumstances.

\begin{remark}
    \label{remark:Reduction}%
    It is clear that the question of $\lie{g}$-covariant Morita theory
    ultimately should result in an understanding of the behaviour of
    Morita theory under \emph{reduction} of star products. Concerning
    the reduction aspect, one has by now a rather good understanding,
    starting from the BRST approach in
    \cite{bordemann.herbig.waldmann:2000a}, see also
    \cite{bordemann:2005a, cattaneo.felder:2007a}. In
    \cite{gutt.waldmann:2010a} the representation theory of the
    reduced algebras was studied in detail, including some aspects of
    Morita theory. The usage of invariant star products (and even
    better: invariant star products with a quantum momentum map) will
    hopefully allow also to treat the Morita theory of star products
    on singular quotients, see
    e.g.\cite{bordemann.herbig.pflaum:2007a, pflaum:2003a}.
\end{remark}

% 
% references
%

\begin{footnotesize}
    \renewcommand{\arraystretch}{0.5} 
%    \bibliographystyle{ewde}
%    \bibliography{dqarticle,dqbook,dqprocentry,dqproceeding,preprints,misc,dqthesis,notes,script}

\begin{thebibliography}{10}

\bibitem {ara:1999a}
{\sc Ara, P.: }\newblock {\em {M}orita equivalence for rings with involution}.
\newblock Alg. Rep. Theo.  {\bf 2} (1999), 227--247.

\bibitem {bayen.et.al:1978a}
{\sc Bayen, F., Flato, M., Fr{{\o}}nsdal, C., Lichnerowicz, A., Sternheimer,
  D.: }\newblock {\em Deformation Theory and Quantization}.
\newblock Ann. Phys.  {\bf 111} (1978), 61--151.

\bibitem {benabou:1967a}
{\sc B{\'e}nabou, J.: }\newblock {\em Introduction to Bicategories}.
\newblock In: {\em Reports of the Midwest Category Seminar},   1--77.
  Springer-Verlag, 1967.

\bibitem {bertelson.bieliavsky.gutt:1998a}
{\sc Bertelson, M., Bieliavsky, P., Gutt, S.: }\newblock {\em Parametrizing
  Equivalence Classes of Invariant Star Products}.
\newblock Lett. Math. Phys.  {\bf 46} (1998), 339--345.

\bibitem {bordemann:2005a}
{\sc Bordemann, M.: }\newblock {\em (Bi)Modules, morphisms, and reduction of
  star-products: the symplectic case, foliations, and obstructions}.
\newblock Trav. Math.  {\bf 16} (2005), 9--40.

\bibitem {bordemann.herbig.pflaum:2007a}
{\sc Bordemann, M., Herbig, H.-C., Pflaum, M.~J.: }\newblock {\em A homological
  approach to singular reduction in deformation quantization}.
\newblock In: {\sc Ch{\'e}niot, D., Dutertre, N., Murolo, C., Trotman, D.,
  Pichon, A. (eds.): }\newblock {\em Singularity theory},   443--461. World
  Scientific Publishing, Hackensack, 2007.
\newblock Proceedings of the {S}ingularity {S}chool and {C}onference held in
  {M}arseille, {J}anuary 24--{F}ebruary 25, 2005. Dedicated to Jean-Paul
  Brasselet on his 60th birthday.

\bibitem {bordemann.herbig.waldmann:2000a}
{\sc Bordemann, M., Herbig, H.-C., Waldmann, S.: }\newblock {\em BRST
  Cohomology and Phase Space Reduction in Deformation Quantization}.
\newblock Commun. Math. Phys.  {\bf 210} (2000), 107--144.

\bibitem {bursztyn:2002a}
{\sc Bursztyn, H.: }\newblock {\em Semiclassical geometry of quantum line
  bundles and {M}orita equivalence of star products}.
\newblock Int. Math. Res. Not.  {\bf 2002}.16 (2002), 821--846.

\bibitem {bursztyn.dolgushev.waldmann:2009a:pre}
{\sc Bursztyn, H., Dolgushev, V., Waldmann, S.: }\newblock {\em Morita
  equivalence and characteristic classes of star products}.
\newblock Preprint  {\bf arXiv:0909:4259} (September 2009), 57 pages.
\newblock To appear in Crelle's J. reine angew. Math.

\bibitem {bursztyn.waldmann:2000b}
{\sc Bursztyn, H., Waldmann, S.: }\newblock {\em Deformation Quantization of
  Hermitian Vector Bundles}.
\newblock Lett. Math. Phys.  {\bf 53} (2000), 349--365.

\bibitem {bursztyn.waldmann:2001a}
{\sc Bursztyn, H., Waldmann, S.: }\newblock {\em Algebraic Rieffel Induction,
  Formal Morita Equivalence and Applications to Deformation Quantization}.
\newblock J. Geom. Phys.  {\bf 37} (2001), 307--364.

\bibitem {bursztyn.waldmann:2002a}
{\sc Bursztyn, H., Waldmann, S.: }\newblock {\em The characteristic classes of
  {M}orita equivalent star products on symplectic manifolds}.
\newblock Commun. Math. Phys.  {\bf 228} (2002), 103--121.

\bibitem {bursztyn.waldmann:2004a}
{\sc Bursztyn, H., Waldmann, S.: }\newblock {\em Bimodule deformations,
  {P}icard groups and contravariant connections}.
\newblock K-Theory  {\bf 31} (2004), 1--37.

\bibitem {bursztyn.waldmann:2005b}
{\sc Bursztyn, H., Waldmann, S.: }\newblock {\em Completely positive inner
  products and strong {M}orita equivalence}.
\newblock Pacific J. Math.  {\bf 222} (2005), 201--236.

\bibitem {bursztyn.waldmann:2005a}
{\sc Bursztyn, H., Waldmann, S.: }\newblock {\em Hermitian star products are
  completely positive deformations}.
\newblock Lett. Math. Phys.  {\bf 72} (2005), 143--152.

\bibitem {calon:2010a}
{\sc Calon, P.: }\newblock {\em Klassischer {L}imes f{\"u}r {$H$}-kovariante
  {D}arstellungstheorie}.
\newblock master thesis, Fakult{\"{a}}t f{\"{u}}r Mathematik und Physik,
  Physikalisches Institut, Albert-Ludwigs-Universit{\"{a}}t, Freiburg, 2010.

\bibitem {cattaneo.felder:2007a}
{\sc Cattaneo, A.~S., Felder, G.: }\newblock {\em Relative formality theorem
  and quantisation of coisotropic submanifolds}.
\newblock Adv. Math.  {\bf 208} (2007), 521--548.

\bibitem {dolgushev:2005a}
{\sc Dolgushev, V.~A.: }\newblock {\em Covariant and equivariant formality
  theorems}.
\newblock Adv. Math.  {\bf 191} (2005), 147--177.

\bibitem {fedosov:1996a}
{\sc Fedosov, B.~V.: }\newblock {\em Deformation Quantization and Index
  Theory}.
\newblock Akademie Verlag, Berlin, 1996.

\bibitem {gerstenhaber:1964a}
{\sc Gerstenhaber, M.: }\newblock {\em On the Deformation of Rings and
  Algebras}.
\newblock Ann. Math.  {\bf 79} (1964), 59--103.

\bibitem {gutt.rawnsley:1999a}
{\sc Gutt, S., Rawnsley, J.: }\newblock {\em Equivalence of star products on a
  symplectic manifold; an introduction to Deligne's {\v{C}}ech cohomology
  classes}.
\newblock J. Geom. Phys.  {\bf 29} (1999), 347--392.

\bibitem {gutt.waldmann:2010a}
{\sc Gutt, S., Waldmann, S.: }\newblock {\em Involutions and Representations
  for Reduced Quantum Algebras}.
\newblock Adv. Math.  {\bf 224} (2010), 2583--2644.

\bibitem {jansen:2006a}
{\sc Jansen, S.: }\newblock {\em {$H$}-{\"A}quivariante {M}orita-{\"A}quivalenz
  und {D}eformationsquantisierung}.
\newblock PhD thesis, Fakult{\"{a}}t f{\"{u}}r Mathematik und Physik,
  Physikalisches Institut, Albert-Ludwigs-Universit{\"{a}}t, Freiburg, November
  2006.

\bibitem {jansen.neumaier.waldmann:2005a:pre}
{\sc Jansen, S., Neumaier, N., Waldmann, S.: }\newblock {\em Covariant {M}orita
  Equivalence of Star Products}.
\newblock In preparation.

\bibitem {jansen.waldmann:2006a}
{\sc Jansen, S., Waldmann, S.: }\newblock {\em The {$H$}-covariant strong
  Picard groupoid}.
\newblock J. Pure Appl. Alg.  {\bf 205} (2006), 542--598.

\bibitem {jurco.schupp.wess:2002a}
{\sc Jur{\v{c}}o, B., Schupp, P., Wess, J.: }\newblock {\em Noncommutative Line
  Bundles and Morita Equivalence}.
\newblock Lett. Math. Phys.  {\bf 61} (2002), 171--186.

\bibitem {kontsevich:2003a}
{\sc Kontsevich, M.: }\newblock {\em Deformation Quantization of {P}oisson
  manifolds}.
\newblock Lett. Math. Phys.  {\bf 66} (2003), 157--216.

\bibitem {lam:1999a}
{\sc Lam, T.~Y.: }\newblock {\em Lectures on Modules and Rings}, vol. 189 in
  {\em Graduate Texts in Mathematics}.
\newblock Springer-Verlag, Berlin, Heidelberg, New York, 1999.

\bibitem {pflaum:2003a}
{\sc Pflaum, M.~J.: }\newblock {\em On the deformation quantization of
  symplectic orbispaces}.
\newblock Diff. Geom. Appl.  {\bf 19} (2003), 343--368.

\bibitem {rieffel:1974b}
{\sc Rieffel, M.~A.: }\newblock {\em Morita equivalence for {$C^*$}-algebras
  and {$W^*$}-algebras}.
\newblock J. Pure. Appl. Math.  {\bf 5} (1974), 51--96.

\bibitem {schmuedgen:1990a}
{\sc Schm{\"{u}}dgen, K.: }\newblock {\em Unbounded Operator Algebras and
  Representation Theory}, vol.~37 in {\em Operator Theory: Advances and
  Applications}.
\newblock Birkh{\"{a}}user Verlag, Basel, Boston, Berlin, 1990.

\bibitem {waldmann:2002b}
{\sc Waldmann, S.: }\newblock {\em {M}orita equivalence of {F}edosov star
  products and deformed {H}ermitian vector bundles}.
\newblock Lett. Math. Phys.  {\bf 60} (2002), 157--170.

\bibitem {waldmann:2007a}
{\sc Waldmann, S.: }\newblock {\em Poisson-{G}eometrie und
  {D}eformationsquantisierung. {E}ine {E}inf{\"u}hrung}.
\newblock Springer-Verlag, Heidelberg, Berlin, New York, 2007.

\end{thebibliography}

\end{footnotesize}
\end{document}